\documentclass{article}
\usepackage{amsmath}
\usepackage{amssymb}

\title{Truth and meaningfulness}
\author{Nik Weaver }
\date{July 2025}

\begin{document}

\maketitle

\section{Tarski's undefinability theorem}

The general form of Tarski's undefinability theorem states that no consistent formal system that meets a minimal standard of expressiveness can contain a truth predicate that applies to its own sentences. If it did, we would be able to formulate a liar sentence and get a contradiction (\cite{tarski}, p.\ 273).\footnote{This paper is a condensed version of the account of truth advanced in my recent book \cite{weaver2}.}

Here ``a truth predicate that applies to its own sentences'' means a predicate $T$ that verifies Tarski's ``convention T'', the schematic biconditional $$T(\mbox{``$\phi$''}) \,\,\,\,\leftrightarrow\,\,\,\, \phi,$$ for every sentence $\phi$ of the language. As this would include sentences in which $T$ itself appears, the impossibility of having such a predicate should not be surprising. The stated requirement is clearly too strong. Once a potentially circular truth predicate is in play, we cannot take it for granted that every syntactically correct sentence will be meaningful, but convention T can only be meaningfully applied to meaningful sentences. If $\phi$ is meaningless then so is the biconditional $T(\mbox{``$\phi$''}) \leftrightarrow \phi$.

Thus, in order to develop a workable theory of global self-applicative truth, it will be necessary to keep track of which sentences are meaningful, and only demand that convention T apply to them.

Tarski's reasoning also shows that we cannot assume every meaningful sentence has a definite truth value. For suppose there was a global truth value function $\tau$, defined on all meaningful sentences and taking the values $t$ and $f$, such that we could affirm $$(\tau(\mbox{``$\phi$''}) = t) \,\,\,\,\leftrightarrow\,\,\,\, \phi$$ for every meaningful sentence $\phi$. Then in any language that could reference $\tau$ and was sufficiently expressive to admit basic self-referential constructions, we could formulate a sentence that says that $\tau$ takes the value $f$ on itself. This sentence would undeniably be meaningful, and either possible truth value for it would lead to a contradiction. This shows that if every meaningful sentence had a definite truth value, then we could prove a contradiction. The evident conclusion is that we have to use intuitionistic logic to reason about global self-applicative truth.

\section{Axiomatizing meaningfulness}

The notion of meaningfulness satisfies some simple axioms. In the propositional setting, we can affirm that meaningfulness is {\it compositional}: a complex proposition is meaningful if and only if each of its atomic subformulas is meaningful. Writing $M(\phi)$ for ``$\phi$ is meaningful'', this is expressed by the equivalences $$(M(\phi) \wedge M(\psi)) \,\,\,\,\leftrightarrow\,\,\,\, M(\ulcorner\phi \wedge \psi\urcorner) \,\,\,\,\leftrightarrow\,\,\,\, M(\ulcorner\phi \vee \psi\urcorner) \,\,\,\,\leftrightarrow\,\,\,\, M(\ulcorner\phi \to \psi
\urcorner)$$ (i.e., each of these four conditions implies the others). In intuitionistic logic we may define the negation $\neg \phi$ of $\phi$ to be the formula $\phi \to (0 = 1)$; since ``0 = 1'' is certainly meaningful, the equivalence $M(\phi) \leftrightarrow M(\ulcorner\neg \phi\urcorner)$ follows.

Moving to the predicate setting, if the formula $\phi$ contains free variables let us
interpret $M(\phi)$ as ``$\phi$ becomes a meaningful sentence when its free variables are replaced by names for any objects''. We can say that such a formula is {\it universally meaningful}. The propositional axioms stated above remain valid, and quantified formulas satisfy the equivalences $$M(\phi) \,\,\,\,\leftrightarrow\,\,\,\, M(\ulcorner(\forall x)\phi\urcorner) \,\,\,\,\leftrightarrow\,\,\,\, M(\ulcorner(\exists x)\phi\urcorner).$$

Any ascription of meaningfulness is always itself meaningful. Thus we have the final axiom $M(\ulcorner M(\mbox{``$\phi$''})\urcorner)$, with $\phi$ ranging over all formulas.

Universal meaningfulness might seem like an unreasonably strong requirement, as formulas with free variables are typically intended to apply only to a limited range of objects, not to all objects without restriction. But this is less of a problem than one might think, because a relation symbol that is initially defined only for variables lying in some special ranges can be extended by simply decreeing that it has the meaning ``0 = 1'' when any of the variables falls outside its intended range. Using this expedient we can convert formulas that are only meaningful in limited settings into universally meaningful ones. Quantification over the originally intended ranges can be achieved using expressions like $(\forall x: \rho(x))\phi \equiv (\forall x)(\rho(x) \to \phi)$ and $(\exists x: \rho(x))\phi \equiv (\exists x)(\rho(x) \wedge \phi)$.

However, since we are reasoning intuitionistically there is an important caveat here. To see the problem, suppose $r$ is a unary relation symbol that is initially defined only for those values of $x$ which satisfy some condition $\rho(x)$. We want to extend it to a unary relation $\tilde{r}(x)$ that is meaningful for arbitrary values of $x$ by setting $\tilde{r}(x) \leftrightarrow r(x)$ or $\tilde{r}(x) \leftrightarrow 0=1$ depending on whether $\rho(x)$ holds. But in order for this to work we need $\rho$ to satisfy excluded middle in the form $$(\forall x)(\rho(x) \vee \neg \rho(x)).$$ Otherwise we will not have given $\tilde{r}(x)$ a definite meaning for all $x$. If for some value of $x$ we do not know which, if either, of $\rho(x)$ or $\neg\rho(x)$ holds, then we do not know which prong of the definition of $\tilde{r}$ to use and therefore cannot consider $\tilde{r}(x)$ to be a meaningful statement.

In short, a formula that is initially defined only when its variables lie in certain ranges can be straightforwardly extended to be meaningful on all values of its variables --- but only if those initially given ranges are {\it definite} in the sense that every object definitely does or does not belong to each of them.

\section{Axiomatizing assertibility}

The other essential tool we will need is the notion of constructive truth, which I understand as rational compulsion or {\it assertibility}. Write $A(\phi)$ for ``$\phi$ is assertible''. This should be taken to mean that we are rationally compelled to affirm the sentence $\phi$ (or the universal closure of $\phi$, if it contains free variables).

In constructivism we understand that the mathematical universe cannot be fully captured axiomatically. To the contrary, we have the ability --- even, one might say, a rational compulsion --- to go beyond any given partial axiomatization, by affirming its consistency, for example. This gives assertibility an open-ended quality that cannot be captured in a single axiomatization. Nonetheless, we can still identify several axioms that an assertibility predicate should satisfy.

For instance, we have axioms that allow us to reason under $A$ using Hilbert-style deduction. Let ``${\rm Log}(\phi)$'' be the statement that the formula $\phi$ is an instance of one of the eleven schemes in the standard axiomatization of the intuitionistic predicate calculus. The validity of these schemes is expressed by the axiom
\begin{eqnarray*}
    M(\phi) \,\,\wedge\,\, {\rm Log}(\phi) &\to& A(\phi)
    \end{eqnarray*}
    and the three standard deduction rules manifest as the axioms
\begin{eqnarray*}
    A(\phi) \,\,\wedge\,\, A(\ulcorner\phi \to \psi\urcorner) &\to& A(\psi)\\
    A(\ulcorner\phi \to \psi(x)\urcorner) &\to& A(\ulcorner\phi \to (\forall y)\psi(y)\urcorner)\\
    A(\ulcorner\psi(x) \to \phi\urcorner) &\to& A(\ulcorner(\exists y)\psi(y) \to \phi\urcorner)
\end{eqnarray*}
where $y$ can be either $x$ or any variable that does not appear freely in $\psi$.

We include a meaningfulness premise in the first displayed axiom because we also want to affirm an axiom which states that every assertible formula is universally meaningful, $$A(\phi) \,\,\,\,\to\,\,\,\, M(\phi).$$ If $\phi$ is not meaningful then even a syntactic tautology such as $\phi \wedge \psi \to \psi$ is not meaningful and cannot be affirmed, for any $\psi$.

The deduction rule axioms are justified by the fact that the corresponding deduction rules are intuitionistically valid. Whenever we have a proof of $\phi$ and a proof of $\phi \to \psi$, we can combine them and invoke modus ponens to get a proof of $\psi$. Thus if $\phi$ and $\phi \to \psi$ are both assertible, then so is $\psi$. The other two deduction rule axioms have similar justifications.

Along similar lines, we have the ``$\forall$-capture'' axiom $$(\forall x: \rho(x))A(\phi(\hat{x})) \,\,\,\,\to\,\,\,\, A(\ulcorner(\forall x): \rho(x)\phi(x)\urcorner),$$ valid provided the condition $\rho(x)$ describes a set and not a proper class. Any set of proofs of the individual instances $\phi(\hat{x})$ can be amalgamated into a single proof of $(\forall x: 
\rho(x))\phi(x)$. Here $\hat{x}$ is assumed to be a name for the corresponding value of $x$.

Finally, we also have the axiom $M(\ulcorner A(\mbox{``$\phi$''})\urcorner)$; ascriptions of assertibility are always meaningful.

\section{Capture and release}

Assertibility figures in the constructive version of the liar paradox that arises from the {\it assertible liar sentence} $\Lambda_a \equiv$ ``$ A(\Lambda_a) \to (0=1)$''.

Convention T is not constitutive of assertibility in the way that it is constitutive of classical truth. We have to consider which direction, if either, of Tarski's biconditional can be affirmed of $A$: {\it capture}, $$\phi \,\,\,\,\to\,\,\,\, A(\mbox{``$\phi$''}),$$ and/or {\it release}, $$A(\mbox{``$\phi$''}) \,\,\,\,\to\,\,\,\, \phi.$$

The capture scheme, with $\phi$ any meaningful sentence, can be seen as a version of Dummett's {\it principle K}, ``if a statement is true, it must be in principle possible to know that it is true'' \cite{dummett}, and might be defensible on the sort of verificationist grounds he champions. It can also be viewed as a consequence of the constructivist conception, mentioned in the last section, of a mathematical universe that is made increasingly explicit --- more poetically, {\it conjured into being} --- by the adoption of successively stronger axioms. In this picture a statement cannot be true until axioms from which it is provable are in place.

But that is not to say that, conversely, everything we are rationally compelled to accept actually is the case: we do not have the right to affirm the release scheme, not even the special case where $\phi$ is ``0=1''. The idea of there being a rationally compelling proof that zero equals one is certainly outrageous, but ``outrageous'' does not entail ``false''. We have no external guarantee that our ever-expanding set of axioms will remain consistent indefinitely.

We may be firmly convinced that there can be no rationally compelling proof of a falsehood. The problematic circularity arises when this conviction is used to justify the release law and that law is subsequently used in arguments that are then recognized as rationally compelling. The release law is, in this way, circular, and as the constructive liar paradox shows, it is a vicious circularity.

For without the release law, the assertible liar sentence $\Lambda_a$ is no longer paradoxical. From the assumption that $\Lambda_a$ is not assertible we can derive ``$0=1$'', so $\Lambda_a$ is not not assertible, but from the assumption that $\Lambda_a$ is assertible we can only derive that ``$0=1$'' is assertible. Without release we cannot infer from this that $\Lambda_a$ is not assertible. Thus $A(\Lambda_a)$ is {\it anomalous}: it is not not the case, and if it is the case then ``$0=1$'' is assertible. This comes awfully close to being a contradiction, but there is no actual contradiction here.

As long as we admit the possibility that there could be a rationally compelling proof of a falsehood, we cannot adopt the release law, and this blocks the constructive liar paradox. The moment we decide such a proof cannot exist, we are forced to accept the $A(\mbox{``0=1''}) \to 0=1$ instance of the release law and to then acknowledge the constructive liar paradox as a valid proof that $0=1$.

Thus, we can augment the axioms for $A$ given in Section 3 with the capture scheme, but not the release scheme. To the extent that a proof of $A(\phi)$ constructively provides a proof of $\phi$, we can affirm the {\it release rule} ``infer $\phi$ from $A(\phi)$'', but this can only be done in special settings, not as a general principle.

\section{Subjunctive reasoning}

If $\phi$ is a meaningful sentence then we define $T(\phi)$ (``$\phi$ is true'') using convention T. But for arbitrary sentences which might not be meaningful, we must make the definition {\it subjunctively}: if $\phi$ were meaningful, then we could affirm the corresponding instance of convention T. That is, we affim the global statement $$M(\phi) \,\,\,\,\to\,\,\,\, A(\ulcorner T(\mbox{``$\phi$''})\leftrightarrow \phi\urcorner).$$ We can reason subjunctively, under $A$, and derive global facts about truth such as $$(M(\phi) \wedge M(\psi)) \,\,\,\,\to\,\,\,\, A(\ulcorner T(\mbox{``$\phi \wedge \psi$''}) \,\,\leftrightarrow\,\, T(\mbox{``$\phi$''}) \wedge T(\mbox{``$\psi$''})\urcorner),$$ with $\phi$ and $\psi$ ranging over all sentences. For any particular sentences $\phi$ and $\psi$ that are known to be meaningful, we can apply the release rule to get $T(\ulcorner \phi \wedge \psi\urcorner) \leftrightarrow (T(\phi) \wedge T(\psi))$. For any set-sized language all of whose sentences are known to be meaningful, $\forall$-capture followed by the release rule would yield $$(\forall \phi,\psi: {\rm Sent}(\phi) \wedge {\rm Sent}(\psi))(T(\ulcorner \phi \wedge \psi\urcorner) \leftrightarrow (T(\phi) \wedge T(\psi)))$$ where ``${\rm Sent}(\phi)$'' says that $\phi$ is a sentence of the given language. As we discussed in Section 2, the universal quantifier is meaningful here because ${\rm Sent}(\phi)$ is a definite predicate and $T(\phi)$ is meaningful for all $\phi$ in this range.

Also as in Section 2, we can define truth for meaningless sentences by setting
$$\neg M(\phi) \,\,\to\,\, A(\ulcorner \neg T(\mbox{``$\phi$''})\urcorner).$$ But this does not make $T(\phi)$ universally meaningful. There could be sentences that are neither definitely meaningful nor definitely meaningless, and for such a sentence we would not know which prong of the definition of $T$ to use so we could not  affirm that $T(\phi)$ has a well-defined meaning.

In particular, it is not clear whether the {\it liar sentence} $\Lambda \equiv$ ``$\neg T(\Lambda)$'' is definitely either meaningful or meaningless. What we can say is that it is not meaningless, i.e., not not meaningful, and if it were meaningful then we could prove that zero equals one. In the terminology of Section 4, the statement $M(\Lambda)$ is anomalous.

Is the liar sentence true or not? The question is senseless until we have said what it means for a sentence to be true, but we can only do this for sentences that definitely are or are not meaningful. Until we have determined whether the liar sentence is meaningful, we cannot even make sense of the question.

The mode of reasoning outlined here involving assertibility, truth, and meaningfulness yields concrete results, such as that $A(\Lambda_a)$ and $M(\Lambda)$ are both anomalous sentences, and is provably consistent (\cite{weaver2}, Theorems 6.1 and 6.2). The consistency of the axioms for assertibility alone, presented in Sections 3 and 4, was previously shown in (\cite{weaver1}, Theorem 5.1). Indeed, when Peano arithmetic is augmented with a self-applicative assertibility predicate we get a consistent formal system which proves of itself that it is both assertibly consistent --- the statement of its consistency is assertible, and assertibly sound --- every theorem of the system is assertible (\cite{weaver1}, Theorems 5.2 and 5.3).

\section{Holding and concepts}

A predicate $\phi$ {\it holds} of an object $x$ if the sentence $\phi(\hat{x})$ is true. Here we can either postulate that every object $x$ has a canonical name $\hat{x}$, or we can work within a definite range (cf.\ Section 2) in which this is known to be the case.

Schematically, we can define holding using ``convention H'', $$H(\mbox{``$\phi$''}, x) \leftrightarrow \phi(x),$$ valid for all $\phi$ and $x$ such that $\phi(\hat{x})$ is meaningful. The global definition is $$M(\phi(\hat{x})) \to A(\ulcorner H(\mbox{``$\phi$''}, \hat{x}) \leftrightarrow \phi(\hat{x})\urcorner).$$

The {\it Russell predicate} $R(\phi) \equiv$ ``$\neg H(\phi, \phi)$'' is paradoxical in essentially the same way as the liar sentence: the statement $M(R(R))$ that $R$ is meaningful when applied to itself is anomalous. We also have a constructive version of the paradox with $R_a(\phi) \equiv$ ``$\neg A(\ulcorner\phi(\mbox{``$\phi$''})\urcorner)$''; here $R_a(R_a)$ is definitely meaningful (as $\Lambda_a$ is) but $A(R_a(R_a))$ is anomalous.

Now let $\phi(x)$ and $\psi(x)$ be formulas with exactly one free variable and suppose they are universally meaningful in the sense of Section 2. We say that they {\it instantiate the same concept} and write $\phi \sim \psi$ if they hold of the same objects. Thus $$(\phi \sim \psi) \,\,\,\,\leftrightarrow\,\,\,\, T(\ulcorner(\forall x)(\phi(x) \leftrightarrow \psi(x))\urcorner).$$ This condition is schematic, stated separately for each universally meaningful $\phi$ and $\psi$. The general rule, stated for all $\phi$ and $\psi$ simultaneously, is $$(M(\phi) \wedge M(\psi))\,\,\,\,\to\,\,\,\, A(\ulcorner(\mbox{``$\phi$''} \sim \mbox{``$\psi$''}) \,\,\leftrightarrow\,\, T(\ulcorner(\forall x)(\phi(x) \leftrightarrow \psi(x))\urcorner)\urcorner).$$

The object $x$ {\it falls under} the concept instantiated by $\phi$ if $\phi$ holds of $x$.

There is no ``Russell concept'' under these definitions, because the predicate ``is a universally meaningful predicate that does not hold of the concept it instantiates'' is not universally meaningful. The part about ``not holding of the concept it instantiates'' is meaningful for any universally meaningful predicate, but this is not a definite range. We cannot say that every predicate definitely is or is not universally meaningful. (Again, cf.\ Section 2.)
 
\bibliographystyle{amsplain}

\end{document}